\documentclass[12pt]{article}
\usepackage{amssymb}
\usepackage{color}

\newtheorem{definition}{Definition}
\newtheorem{theorem}[definition]{Theorem}
\newtheorem{proposition}[definition]{Proposition}
\newtheorem{example}{Example}

\newcommand\unit{\hbox{\rm 1\kern-2.8truept l}}
\newcommand\re{{\Re}\kern-1pt e}
\newcommand\im{{\Im}\kern-1pt m}

\newcommand{\B}{\mathcal{B}(\mathsf{h})}
\newcommand\Lform{{\mathcal{L}}\kern-7.56pt\raise1.0pt\hbox{$-$}}

\newcommand{\h}{\mathsf{h}}
\newcommand{\T}{\mathcal{T}}
\newcommand{\Ll}{\mathcal{L}}
\newcommand{\mi}{{\mathrm{i}}}

\newcommand{\ketbra}[2]{\left| #1\right\rangle\left\langle #2\right|}
\newcommand{\binomial}[2]{ \left(
\begin{array}{c} {#1}\\{#2} \end{array} \right)} 

\begin{document}

\title{On the relationship between a quantum Markov semigroup  
and its representation via linear stochastic Schr\"odinger equations}   

\author{FRANCO FAGNOLA\footnote{Dipartimento di Matematica, Politecnico di Milano, 
Piazza Leonardo da Vinci 32, I-20133 Milano, Italy
{\tt franco.fagnola@polimi.it}}  \ 
and CARLOS MORA\footnote{CI$^2$MA and Departamento de Ingenier\'{\i}a Matem\'{a}tica, Universidad de Concepci\'on,
Barrio Universitario, Avenida     Esteban Iturra s/n, 4089100 , 
Casilla 160-C 
Concepci\'on, Chile {\tt cmora@ing-mat.udec.cl}} }  

\date{ }          

\maketitle

\abstract{A quantum Markov semigroup can be represented via 
classical diffusion processes solving a stochastic Schr\"odinger
equation. In this paper we first prove that a quantum Markov 
semigroup is irreducible if and only if classical diffusion 
processes are total in the Hilbert space of the system. Then we 
study the relationship between irreducibility of a quantum 
Markov semigroup and properties  of these diffusions such 
as accessibility, the Lie algebra  rank condition, and irreducibility. 
We prove that all these properties are, in general, weaker than 
irreducibility of the quantum Markov semigroup, nevertheless, 
they are equivalent for some important classes of semigroups.}

\bigskip

{\bf Keywords.} Open quantum systems, quantum Markov semigroups, 
stochastic Schr\"odinger equations, irreducibility, support of quantum 
states, control.


\section{Introduction}

A quantum Markov semigroup (QMS) $\T$  is a weakly$^*$-continuous
semigroup $(\T_t)_{t\ge 0}$ of completely positive, identity preserving, 
normal maps on a von Neumann algebra. In this paper, we will only be concerned with QMS on a matrix algebra which are norm-continuous.  

These QMS semigroups were introduced in the seventies (as quantum 
dynamical semigroups) to model the irreversible evolution of an open 
quantum system and are now an important tool to investigate quantum 
systems and quantum stochastic processes (see \cite{BarcGreg,BaLiSk,BrPe,FaWi,GoSi,JPW,LiWi} and the references 
therein). The representation of QMS via solutions of classical 
stochastic differential equations, already noticed by A.V. Skorohod 
\cite{Skoro}, is also well-known and plays a key role in quantum 
trajectory  theory  (see, e.g. \cite{BarcGreg} section 3.2.3, 
\cite{BaPePe,BaBe,Slava,BrPe,Kolok,KuMa} and the references  
therein). These equations, called stochastic Schr\"odinger equations 
(SSE) (see \cite{FaMo,Mora-JAF2008,Mora-AP2013,Mora-RR-IDAQP2007,
Mora-RR-AAA2008,Pellegrini-AP2008} for recent results), 
turn out to be very useful 
to study open quantum systems through the interplay between 
classical and quantum stochastic analysis. 

The aim of this note is to study the relationship between 
irreducibility of a QMS and diffusion processes solving the 
associated SSE driven by independent Brownian motions. 
Our motivation is to establish a bridge between classical and 
quantum stochastics, however, these results may turn out to be 
useful in the study of open quantum systems and their numerical 
simulations via SSEs since irreducibility enables one to apply powerful 
results from ergodic theory. Moreover, we want to find the range 
of solutions to SSE because it can be thought of as the set of 
reachable (random) states in a continuous measurement.

We first give a new characterisation 
(Theorem \ref{th:irred-iff-Sxi=h}) irreducible QMS 
by a multiple commutator condition looking like the celebrated 
Lie algebra rank condition (LARC) and H\"ormander condition.

Then we prove 
our main result (Theorem \ref{th:QMS-irred-iff-SSE-tot}); 
a QMS is irreducible if and only if the associated  diffusion processes 
via SSEs are total in the Hilbert space of the system. 

Moreover, we study the relationship with other properties such as 
accessibility, the Lie algebra rank condition, and irreducibility.  
We prove that irreducibility of a QMS is, in general, a weaker 
property than irreducibility of diffusions solving the associated SSEs. 
It is also weaker of the LARC and H\"ormander condition,  
although equivalent for some important classes of semigroups. 
We thus find a quantum version of these classical conditions.

The paper is organized as follows. In section \ref{sect:irred-QMS} 
we present a short account of the main results on irreducible QMSs 
and describe the support projection at time $t$ of a state evolving 
under the action of a QMS  together  with the characterisation based 
on the multiple commutator condition 
(Theorem \ref{th:irred-iff-Sxi=h}). 
In section \ref{sect:SSE} we introduce SSE, driven by independent 
Brownian motions, and prove (Theorem \ref{th:QMS-irred-iff-SSE-tot}) that a QMS is irreducible if and only if associated diffusion processes 
are total in the Hilbert space of the system. Then we turn our 
attention to the range of  diffusion processes showing by simple 
argument and examples (Example \ref{ex:Lanti-sa}) that  we 
cannot expect these diffusion processes to be irreducible.
In section \ref{sect:range-sol-SSE} we discuss the Stroock and 
Varadhan support theorem and the LARC condition showing 
(Theorem \ref{th:LARC} and counterexample \ref{ex:3x3}) 
that our multiple commutator condition, equivalent to 
irreducibility of a QMS, is indeed weaker that the LARC 
condition. Finally, in section \ref{sect:Generic-QMSs}, we 
show that both the LARC condition and the multiple 
commutator condition hold for generic QMSs.

\section{Irreducible QMS}\label{sect:irred-QMS}

Let $\h=\mathbb{C}^d$ and let $\T$ be the QMS on 
the algebra $M_d(\mathbb{C})$ of $d\times d$ matrices 
generated by 
\begin{eqnarray}\label{eq:GKSL}
\Ll(x) & = &  \mi [H,x]
+ \frac{1}{2} \sum_{\ell=1}^m\left( 
- L^*_\ell L_\ell x  +2 L^*_\ell x L_\ell - x L^*_\ell L_\ell\right) 
\nonumber\\
&= & G^*x + \sum_{\ell=1}^m L^*_\ell x L_\ell + x G 
\end{eqnarray}
where $L_1,\dots, L_m, H\in M_d(\mathbb{C})$ with $H$ 
self-adjoint and 
\[
G=-\frac{1}{2}\sum_{\ell=1}^m L_\ell^* L_\ell -\mi H.
\]

The representation (\ref{eq:GKSL}) of the generator $\Ll$ 
is called a Gorini-Kossa\-kow\-ski-Sudarshan-Lindblad (GKSL) 
representation of the generator. 
It is well-known to exist but it is not unique 
 (see \cite{Partha} Theorem 30.16 p. 271).  In particular,
one can always change a GKSL representation by translating 
the operators $L_\ell$  by multiples of the identity operator 
or increasing  $m$ and adding operators $L_{j}$ which are 
multiples of the identity operator.

A representation with the smallest number of operators $L_\ell$, i.e. the minimum $m$,
is called \emph{minimal}. 
In a minimal GKSL representation of $\Ll$, matrices 
$\unit,L_1,\dots,L_m$ are linearly independent 
 (\cite{Partha} Theorem 30.16 p. 271).

\begin{definition}\label{def:S-xi-manifold}
For each non-zero $\xi\in \h$ we denote by $\mathcal{S}(\xi)$ 
the linear span of all vectors of the form
\begin{equation}\label{eq:supp-state-diff}
\xi, \,
\delta_G^{n_1}(L_{\ell_1})\delta_G^{n_2}(L_{\ell_2}) 
\cdots \delta_G^{n_k}(L_{\ell_{k}})\xi ,
\end{equation}
where
$k\ge 1$, $n_1,\dots,n_k\ge 0$ and $ 1 \leq \ell_1,
\dots,\ell_k \leq m$ and 
$\delta_G^{n}$ is defined recursively by $\delta_G^{0}(A) = A$, 
$\delta_G^{n+1}(A) = \left[ G, \delta_G^{n}(A)\right] $.
\end{definition}

The following results are proved in \cite{FFRR-ep} 
Theorem 6 and 7 (see also their extensions in 
\cite{Hachicha}).

\begin{theorem}\label{th:supp-state}
Let $(\T_t)_{t\ge 0}$ be a norm continuous QMS on $\B$ 
with generator $\Ll$  as in (\ref{eq:GKSL}) and let $P_t = 
\hbox{\rm e}^{tG}$. For all unit vector 
$\xi\in\h$ and all $t\ge 0$, the support projection of the state 
$\T_{*t}(\ketbra{\xi}{\xi})$ is the closed linear span of $P_t \xi$ 
and vectors 
\begin{equation}\label{eq:supp-state}
P_{s_1}L_{\ell_1}P_{s_{2}-s_{1}}L_{\ell_2}P_{s_3-s_2} \dots 
P_{s_{n}-s_{n-1}}L_{\ell_n}P_{t-s_n}\xi
\end{equation}
for all $n\ge 1$, $0\le s_1 \le s_{2}\le \dots \le s_{n}\le t$ 
and $\ell_1,\dots,\ell_n\ge 1$. 
\end{theorem}

A simple argument based on the analyticity of the semigroup 
$(P_t)_{t\ge 0}$ (see \cite{FFRR-ep} Theorem 7) leads to 
the following simpler characterisation of the support of 
$\T_{*t}(\ketbra{\xi}{\xi})$ as the linear manifold 
$P_t \,\mathcal{S}(\xi)$.

\begin{theorem}\label{th:supp-state-differential}
Let $(\T_t)_{t\ge 0}$ be a norm continuous QMS on 
$\B$ with generator $\Ll$  as in 
(\ref{eq:GKSL}) and let $P_t = \hbox{\rm e}^{tG}$. 
For all unit vector $\xi\in\h$ and all $t > 0$, the support 
projection of the state $\T_{*t}(\ketbra{\xi}{\xi})$ is the 
linear manifold $P_t \,\mathcal{S}(\xi)$ where $\mathcal{S}(\xi)$ 
is  the linear span of vectors (\ref{eq:supp-state-diff}).
\end{theorem}

\begin{definition}\label{def:irred}
A QMS $\T$ is \emph{irreducible} if there exists no non-trivial 
subharmonic projection $p$ ($\T_t(p)\ge p$ for all $t\ge 0$).
\end{definition}

In an equivalent way, a QMS $\T$ is irreducible there exists 
no non-trivial common invariant subspace for the operators 
$G$ and $L_\ell$ (\cite{FFRR-jmp02} Theorem III.1).

Let $\mathcal{F}(\T)$ be the vector space of fixed points of $\T$
\[
\mathcal{F}(\T) = \left\{ \, x \mid \T_t(x)=x,\, \forall t\ge 0\,\right\}.
\]
It is well-known that, if $\T$ has a faithful invariant state, then 
$\mathcal{F}(\T)$ is a sub-$^*$-algebra of $M_d(\mathbb{C})$. 

Let $\mathcal{N}(\T)$ be the decoherence free algebra 
of $\T$ 
\[
\mathcal{N}(\T) = \left\{ \, x \mid
\T_t(x^*x)=\T_t(x^*)\T_t(x),\, 
\T_t(x x^*)=\T_t(x)\T_t(x^*)\  
\forall t\ge 0\,\right\}.
\]
We refer to \cite{FFRR-idaqp08} for properties of $\mathcal{N}(\T)$.
Both $\mathcal{F}(\T)$ and $\mathcal{N}(\T)$ contain scalar 
multiples of the identity matrix $\unit$; we say that they are trivial 
if they do not contain other matrices, i.e. they coincide with 
$\mathbb{C}\unit$.

We recall the following result on irreducible QMSs 

\begin{theorem}\label{th:irred-FT-trivial}
An irreducible QMS $\T$ on $M_d(\mathbb{C})$ admits a unique 
faithful invariant state. Its fixed point set $\mathcal{F}(\T)$ and  
decoherence free subalgebra $\mathcal{N}(\T)$ are trivial.
\end{theorem}

\noindent{\bf Proof.} By finite-dimensionality, the QMS $\T$ 
admits an invariant state $\rho$ and its support projection is 
subharmonic (see e.g. \cite{FFRR-jmp02} Theorem II.1) and 
non-zero. Thus it must coincide with $\unit$ because $\T$ 
is irreducible and so $\rho$ is faithful.

As a well-known consequence,  $\mathcal{F}(\T)$ is a 
$^*$-subalgebra of $M_d(\mathbb{C})$ because, for 
any $x\in\mathcal{F}(\T)$, by complete positivity, we 
have $\T_t(x^*x)\ge \T_t(x^*)\T_t(x)=x^*x$ and, by the 
invariance of $\rho$, we have 
\[
\hbox{\rm tr}\left(\rho\left( \T_t(x^*x)- \T_t(x^*)\T_t(x)\right) \right)
=\hbox{\rm tr}\left(\rho\left( \T_t(x^*x)- x^*x\right) \right)=0
\] 
Thus $x \in\mathcal{N}(\T)$. 
Since the algebra $\mathcal{N}(\T)$ is trivial by \cite{DFSU} Proposition 14, because $\T$ is irreducible, also $\mathcal{F}(\T)$ 
is trivial.

We finally show that $\rho$ is the unique invariant state 
of $\T$. Indeed, if it is not, then the dimension of the kernel 
of $\Ll_*$ is at  least $2$ and so, since ker$(\Ll_*)$ is the 
orthogonal space of  the range of $\Ll$, it follows that 
the dimension of $R(\Ll)$ is not bigger than $d^2-2$.
This implies that the dimension of the kernel of $\Ll$ is at 
least $2$ contradicting the triviality of $\mathcal{F}(\T)$.
\hfill $\square$

\medskip

The following new characterisation of irreducible QMS can  
be regarded as the starting point of our analysis.

\begin{theorem} \label{th:irred-iff-Sxi=h}
The following are equivalent:
\begin{enumerate}
\item[(1)] the QMS $\T$ is irreducible,
\item[(2)] $\mathcal{S}(\xi)=\h$ for all non-zero $\xi\in \h$,
\end{enumerate}
\end{theorem}

\noindent{\it Proof.} (2) $\Rightarrow$ (1). If the  QMS $\T$ is 
not irreducible there exists a nontrivial subharmonic projection $p$. 
The subspace determined by $p$ is invariant under $G$ and all 
the $L_\ell$ by Theorem III.1 of \cite{FFRR-jmp02}. Therefore for 
all non-zero $\xi$ in the range of $p$, $\mathcal{S}(\xi)$ is contained 
in the range of $p$.  

(1) $\Rightarrow$ (2). If the  QMS $\T$ is irreducible it 
admits a unique faithful invariant state $\rho$ by Theorem 
\ref{th:irred-FT-trivial}. Moreover, since $\mathcal{F}(\T)= \mathcal{N}(\T)$ (indeed both are trivial),  
by  result due to Frigerio and Verri (see Theorem 3.3 
of \cite{FrVe}) for any unit vector $\xi$ in $\mathbb{C}^d$ 
\[
\lim_{t\to\infty} \T_{*t}(\ketbra{\xi}{\xi})=\rho.
\]
By finite dimensionality, it follows that the state
$\T_{*t}(\ketbra{\xi}{\xi})$ is faithful for all $t$ 
bigger than some $t_0<+\infty$.
Hence $P_t S(\xi)=\mathbb{C}^d$  for all $t>t_0$ 
by Theorem \ref{th:supp-state-differential},
and so $S(\xi)=\mathbb{C}^d$ by the 
invertibility of $P_t$. 
\hfill $\square$

\medskip
The previous result yields an algebraic condition that implies 
some qualitative property of a QMS (see \cite{FFRR-2006, FrVe, 
Spohn} for related algebraic conditions implying other properties).

\section{Stochastic Schr\"odinger equations}
\label{sect:SSE}

A linear SSE for the QMS generated by (\ref{eq:GKSL}) 
is the stochastic differential equation 

\begin{equation}\label{eq:SSE}
dX_t(\xi) = GX_t(\xi) dt 
+ \sum_{ \ell = 1}^{m} L_\ell X_t(\xi) dW^{\ell}_t,  
\qquad X_0(\xi)=\xi
\end{equation}
where $\xi \in \mathbb{C}^d$ and $W^{1}, \ldots, W^{m}$ 
are independent real-valued independent Wiener processes 
on a filtered complete probability space 
$(\Omega,\mathfrak{F}, (\mathfrak{F}_t)_{t\ge 0},\mathbb{P})$.

It is well-known that (see e.g. \cite{BarcGreg} Theorem 2.11 p.29)
\[
\langle \eta,\T_t(a) \xi  \rangle = \mathbb{E}
\left[\left\langle X_t(\eta), a X_t(\xi)\right\rangle\right], 
\qquad 
\T_{*t}\left( \left|\xi\right\rangle\left\langle \eta\right|\right)
= \mathbb{E}\left[  \left|X_t(\xi)\right\rangle\left\langle X_t(\eta)\right|\right].
\]
Moreover, since the operators $G,L_\ell$ are bounded and 
$\Ll(\unit)=0$, by \cite{BarcGreg} Theorem 2.11, we have also
\begin{equation}\label{eq:norm-pres}
\mathbb{E}\left[\,  \left\Vert X_t(\xi)\right\Vert^2\, \right]
= \left\Vert \xi\right\Vert^2.
\end{equation}

The initial condition $\xi$ will be always assumed to be non-zero.
In this way we associate with a generator $\Ll$ a diffusion process 
on $\mathbb{C}^d-\{0\}$. In order to investigate the relationship 
between irreducibility the QMS generated by $\Ll$ and the 
diffusion process (\ref{eq:SSE}) we start with 
the following result.

\begin{proposition}\label{prop:chaos-exp}
The random variable $X_t(\xi)$ 
admits the chaos expansion
\begin{eqnarray}
& & \kern-16truept X_t(\xi) = P_t\,\xi  \\
& & \kern-16truept + \kern-8truept
\sum_{n\ge 1,\ell_1,\dots,\ell_n\ge 1} 
\int_0^t dW^{\ell_1}_{s_1} ..
\int_{0}^{s_{n-1}}  \kern-14truept dW^{\ell_n}_{s_n}
P_{t-s_1} L_{\ell_1} P_{s_1-s_2} \cdots 
P_{s_{n\kern-1truept - \kern-1truept1}-s_n}
L_{\ell_n} P_{s_n}\xi. \nonumber
\end{eqnarray}
\end{proposition}

\noindent{\bf Proof.}  Recall that 
 $P_t = \hbox{\rm e}^{  t G } $, i.e., $(P_t)_{t\ge 0}$ is 
 the contraction semigroup generated by $G$. 
For all $t>0$ and $s\in]0,t[$ we 
have
\[
d P_{t-s}X_s(\xi) = 
\sum_{\ell\ge 1} P_{t-s} L_\ell X_s(\xi) dW^{\ell}_s
\]
so that, integrating on $[0,t]$, 
\[
X_t(\xi) = P_t\,\xi + \sum_{\ell\ge 1} \int_0^t P_{t-s} L_\ell X_s(\xi) dW^{\ell}_s.
\]
Iterating this formula $n$ times we can write $X_t(\xi)$ as 
the sum of $P_t\xi$ plus
\begin{eqnarray*}
\sum_{k=1}^n
\sum_{\ell_1,\dots,\ell_k\ge 1} 
\int_0^t dW^{\ell_1}_{s_1} \cdots 
\int_{0}^{s_{k-1}}  \kern-8truept dW^{\ell_k}_{s_k}
P_{t-s_1} L_{\ell_1}\cdots 
P_{s_{k-1}-s_k}L_{\ell_k} P_{s_k}\xi
\end{eqnarray*}
and a remainder $R_n(\xi)$ given by
\[
\sum_{\ell_1,\dots,\ell_{n+1}\ge 1} 
\int_0^t dW^{\ell_1}_{s_1} \cdots 
\int_{0}^{s_{n}}  \kern-8truept dW^{\ell_n}_{s_{n+1}}
P_{t-s_1} L_{\ell_1}\cdots 
P_{s_{n}-s_{n+1}}L_{\ell_{n+1}} X_{s_{n+1}}(\xi).
\]
Therefore, putting 
\[
c:=\max_{1\le \ell \le m}\left\Vert L_\ell\right\Vert ,
\]
from we can write 
$\mathbb{E}\left[\,  \left\Vert R_n(\xi)\right\Vert^2\, \right]$ as 
\begin{eqnarray*}
&  & \kern-12truept\sum_{\ell_1,\dots,\ell_{n+1}\ge 1} 
\int_0^t d{s_1} \cdots 
\int_{0}^{s_{n}} \kern -8truept d{s_{n+1}}\,
\mathbb{E}\left[\, 
 \left\Vert P_{t-s_1} L_{\ell_1}\cdots 
L_{\ell_{n+1}} X_{s_{n+1}}(\xi)\right\Vert^2 \,\right] \\
&  & \kern-12truept\le c^{n+1} m^{n+1}
\sum_{\ell_1,\dots,\ell_{n+1}\ge 1} 
\int_0^t d{s_1} \cdots 
\int_{0}^{s_{n}} \kern -8truept d{s_{n+1}} \\
&  & \kern-8truept = \frac{(c\,m)^{n+1}\, t^{n+1}}{(n+1)!}.
\end{eqnarray*}
The conclusion follows letting $n$ go to infinity.
\hfill$\square$

\medskip

The solution $X_t(\xi)$ to (\ref{eq:SSE}), 
for $t>0$ and $\xi\in\mathbb{C}^d$ fixed, defines a family 
of random vectors on $\Omega$.
We recall that the essential range of a 
$\mathbb{C}^d$-valued random variable $Y$ is the set of all 
$u\in\mathbb{C}^d$ such that 
$\mathbb{P}\{ Y\in \mathcal{U}\} > 0$ 
for each neighbourhood $\mathcal{U}$ of $u$.

\begin{theorem}\label{th:QMS-irred-iff-SSE-tot}
Let $(X_t(\xi)_{t\ge 0}$ be the unique solution to (\ref{eq:SSE}).  
The following are equivalent:
\begin{enumerate}
\item $S(\psi)=\mathbb{C}^d$ for all $\psi\in\mathbb{C}^d
-\{0\}$,
\item for all $t>0$ and all $\xi\in\mathbb{C}^d-\{0\}$ 
the essential range of $X_t(\xi)$ is  total 
in $\mathbb{C}^d$.
\end{enumerate}
\end{theorem}

\noindent{\bf Proof.} 
By the chaos expansion of Proposition \ref{prop:chaos-exp}, for 
all $v\in \mathbb{C}^d$ we have 
\begin{eqnarray*}
& &  \mathbb{E}\left[\,  \left|  \left\langle 
v, X_t(\xi)\right\rangle\right|^2 \,\right] 
= \left\Vert P_t\,\xi\right\Vert^2 \\
& &  + \sum_{n\ge 1,\ell_1,\dots,\ell_n\ge 1} 
\int_0^t d{s_1} \dots
\int_{0}^{s_{n-1}}  \kern-8truept d{s_n}
\left|\left\langle v,
P_{t-s_1} L_{\ell_1} P_{s_1-s_2} \cdots 
L_{\ell_n} P_{s_n}\xi \right\rangle\right|^2.
\end{eqnarray*}
It follows that $\mathbb{P}\left\{\left\langle 
v, X_t(\xi)\right\rangle \not=0 \right\}=0$ if 
and only if $v$ is orthogonal to all vectors $P_t\xi$, 
$P_{t-s_1} L_{\ell_1} P_{s_1-s_2} \cdots 
L_{\ell_n} P_{s_n}\xi$ with $1\le\ell_1,\dots,\ell_n\le m$, 
$0\le s_n \dots \le s_1\le t$.

Now, if 2 holds, $\mathbb{P}\left\{\left\langle 
v, X_t(\xi)\right\rangle \not=0 \right\}>0$ for all 
$v\in\mathbb{C}^d-\{0\}$  and $v$ is not orthogonal to 
all vectors
 $P_t\xi$, $P_{t-s_1} L_{\ell_1} P_{s_1-s_2} \cdots 
L_{\ell_n} P_{s_n}\xi$ with $1\le\ell_1,\dots,\ell_n\le m$, 
$0\le s_n \dots \le s_1\le t$. Thus no non-zero vector is 
orthogonal to $P_t S(\xi)$ by Theorem  
\ref{th:supp-state-differential}. Since $P_t$ is invertible, and 
$S(\xi)$ is a subspace of $\mathbb{C}^d$, it turns out that $S(\xi)=\mathbb{C}^d$.

Conversely, if 1 holds, then, for all non-zero $v\in\mathbb{C}^d$, 
the expectation of the random variable 
$\left|  \left\langle v, X_t(\xi)\right\rangle\right|^2 $ is strictly 
positive by Theorems \ref{th:supp-state}  and 
\ref{th:supp-state-differential} so that 
$X_t(\xi)$ is not orthogonal to $v$ on an event of strictly 
positive probability.
\hfill $\square$

\bigskip 
We can now proceed to study the relationship between irreducibility 
of a QMS and diffusion process $(X_t)_{t\ge 0}$ solving the 
associated SSE first recalling the usual definition. 

\begin{definition}\label{def:irred-diffusion}
The diffusion process $(X_t)_{t\ge 0}$ is called \emph{irreducible} 
on $\mathbb{C}^d-\{0\}$ if, for all $X_0=\xi\in\mathbb{C}^d-\{0\}$  
and all open set $O\subseteq \mathbb{C}^d-\{0\}$, there 
exists $t>0$ such that 
\[
\mathbb{P}\left\{ X_t(\xi) \in O\right\} >0.
\]
\end{definition}

Clearly, even if the QMS associated with $G,L_\ell$ is irreducible,
the diffusion process 
$(X_t(\xi))_{t\ge 0}$ in $\mathbb{C}^d-\{0\}$ may not be. 
This is the case, for instance, when the non-zero vector $\xi$ 
has \emph{real} components and matrices $G,L_\ell$ have 
\emph{real} entries.

This is not just a matter of phase and length of vectors 
$X_t(\xi)$ because, if $d>2$, the diffusion process takes 
values in a manifold of \emph{real} dimension $d$ which is 
strictly smaller than the real dimension $2(d-1)$ of the 
complex projective space $\mathbb{CP}^{d-1}$ obtained 
on taking the quotient with respect to a complex scalar. 

Moreover, these situations indicate that irreducibility of 
the QMS  associated with $G,L_\ell$ is much weaker than 
irreducibility of
\begin{enumerate}  
\item the diffusion process solving (\ref{eq:SSE}) 
in $\mathbb{C}^d-\{0\}$, 
\item the diffusion process on the unit 
sphere  of $\mathbb{C}^d$ obtained by normalization of 
vectors $X_t(\xi)$, 
\item the diffusion process on  the complex 
projective space $\mathbb{CP}^{d-1}$ associated with 
the solution of (\ref{eq:SSE}).
\end{enumerate}

We finish this section by showing another example illustrating 
that irreducibility of (\ref{eq:SSE}) is  stronger 
than irreducibility QMS.

\begin{example}\label{ex:Lanti-sa}
{\rm
Suppose that the operators $L_\ell$ are anti-selfadjoint, namely
$L_\ell^*=-L_\ell$ for $\ell=1,\dots,m$. Then, by the Ito 
formula,
\begin{eqnarray*}
d\left\langle X_t(\xi),X_t(\xi)\right\rangle
& = & \Big\langle X_t(\xi), (G^*+\sum_\ell L_\ell^*L_\ell + G)X_t(\xi)\Big\rangle dt  \\
& + & \sum_{\ell}\left\langle X_t(\xi),(L_\ell^*+L_\ell)X_t(\xi)\right\rangle dW^\ell_t = 0.
\end{eqnarray*}
It follows that $\Vert X_t(\xi)\Vert = \Vert \xi \Vert$ for all 
$t >0$, namely the range of the diffusion process is contained 
in the unit sphere of $\mathbb{C}^d$.

It is not hard, however, to produce an  irreducible QMS with 
anti-selfadjoint operators $L_\ell$. We may consider, for instance
$d=2, m=1$, $L=\mi\sigma_2$, $H=\sigma_3$ where 
$\sigma_1,\sigma_2,\sigma_3$ are the Pauli matrices
\[
\sigma_1 = \left[  \begin{array}{cc} 0 & 1 \\ 1 & 0
                           \end{array}\right], \quad
\sigma_2 = \left[  \begin{array}{cc} 0 & -\mi \\ \mi & 0
                           \end{array}\right], \quad
\sigma_3 = \left[  \begin{array}{cc} 1 & 0 \\ 0 & -1
                           \end{array}\right].
\]
In this case we have 
\[
\delta_{G}(L)=\delta_{-\mi H}(L)=[-\mi\sigma_3,\sigma_2]
=2\sigma_1, \quad \delta_{G}(L)L=2\mi\sigma_3
\] 
so that $S(\xi)=\mathbb{C}^2$ for all $\xi\in\mathbb{C}^2-\{0\}$.
}
\end{example}

Further examples will be discussed in the next section.

\section{The range of solutions to linear SSE}
\label{sect:range-sol-SSE}

Useful tools from control theory are available to study the range 
of solutions to SSE (see the survey \cite{ReyB} Sect.6)
For equation (\ref{eq:SSE}) let us replace the Wiener processes 
$W^\ell$ by piecewise polygonal approximations
\[
W^{\ell,n}_t 
= W^{\ell,n}_{k/n} + (nt - k)
\left( W^{\ell,n}_{(k+1)/n}-W^{\ell,n}_{k/n}\right), \qquad 
\frac{k}{n}\le t \le \frac{k+1}{n}.
\]
A celebrated result by Stroock and Varadhan \cite{StrVar} 
shows that, for any $d\times d$ matrix $\widetilde{G}$, 
the solutions of
\begin{equation}\label{eq:SSE-piecewise}
dX^{(n)}_t(\xi) =  \widetilde{G} X^{(n)}_t(\xi) dt 
+ \sum_{\ell=1}^m L_\ell X^{(n)}_t(\xi)   dW^{\ell,n}_t,  
\qquad X^{(n)}_0(\xi)=\xi 
\end{equation}
converge almost surely to $X_t(\xi)$ uniformly in $t$ on any 
compact interval to the solution of
\begin{equation}\label{eq:SSE-Stratonovich}
dX_t(\xi) = \widetilde{G}X_t(\xi) dt 
+ \sum_{\ell  =1}^m L_\ell X_t(\xi) \circ dW^{\ell}_t,  
\qquad X_0(\xi)=\xi
\end{equation}
where $\circ$ denotes the Stratonovich integral, namely,  
in terms of the Ito integral, to the solution of 
\[
dX_t(\xi) = 
\left(\widetilde{G} +\frac{1}{2}\sum_{\ell=1}^m L_\ell^2
\right) X_t(\xi) dt 
+ \sum_{\ell =1}^m L_\ell X_t(\xi) dW^{\ell}_t,  
\quad X_0(\xi)=\xi.
\]
Thus, choosing
\begin{equation}\label{eq:widetildeG}
\widetilde{G} = G - \frac{1}{2}\sum_{\ell=1}^m L_\ell^2
\end{equation} 
solutions of (\ref{eq:SSE-piecewise}) converge to solutions 
of (\ref{eq:SSE}).

Equation (\ref{eq:SSE-Stratonovich}) has the form
\begin{equation}\label{eq:bilinear-control}
dx_t =  \widetilde{G} x_t  dt 
+ \sum_{\ell=1}^m L_\ell x_t  u_\ell(t)  dt, \qquad x_0=\xi  
\end{equation}
where $u_\ell$ are piecewise constant functions. 
This is an ordinary (non-autonomous) differential equation, 
functions $u_\ell$  are controls and (\ref{eq:bilinear-control} 
is a bilinear control system.

Unfortunately there are no general necessary and sufficient 
condition for deciding when a bilinear system is controllable.

A well-known necessary condition for controllability 
(\cite{Sachkov} Theorem 2.3) is the \emph{Lie algebra rank condition} 

\begin{definition}
The Lie algebra rank condition (LARC) holds if the linear manifold generated 
by vectors 
\begin{equation}\label{eq:subsetLARC}
\widetilde{G}\xi,\,L_\ell\xi,[\widetilde{G},L_\ell]\xi,\,
[L_{\ell_1},L_{\ell_2}]\xi,\,
[\widetilde{G},[\widetilde{G},L_\ell]]\xi,\,
[\widetilde{G},[L_{\ell_1},L_{\ell_2}]]\xi, ...
\end{equation} 
is $\mathbb{C}^d$ for all non zero $\xi$.
\end{definition}

The LARC implies that the control system 
(\ref{eq:bilinear-control}) is \emph{accessible} 
namely the set of points $(x_s)_{0\le s\le t}$ reachable 
from $\xi$ with some choice of piecewise constant 
controls $u_\ell$ contains a non-empty open set 
in $\mathbb{C}^d$ for all $t>0$. 

It is worth noticing here that the linear manifold spanned 
by vectors (\ref{eq:subsetLARC}) may depend on the 
particular choice of the operators $G$ and $L_\ell$ in the 
GKSL representation of the generator $\Ll$ as shows the 
next example. For this reason, from now on, we consider 
only \emph{minimal} GKSL representations of $\Ll$. 

\begin{example}\label{ex:pureHam}
{\rm 
Let $\h = \mathbb{C}^2$, let $H$ be a self-adjoint matrix 
which is not a multiple of the identity $\unit$, $L_\ell=0$ for all 
$\ell\ge 1$ and let $\T$ the QMS on $\mathbb{C}^2$ defined by 
$\T_t(x) = \hbox{\rm e}^{\mi t H} x \hbox{\rm e}^{-\mi t H}$.

The LARC does not hold because the dimension of the 
linear manifold (\ref{eq:subsetLARC}) is at most 1 for 
all non zero $\xi\in\mathbb{C}^2$. However, if we 
consider a GKSL representation of the generator 
$\Ll(x)= \mi [H,x]$ with operators $L_\ell=0$ for 
all $\ell>1$,
\[
L_1 = z \unit, \qquad G = -\frac{|z|^2}{2}\unit -\mi H
\]
for some non zero complex number $z$, the linear manifold 
(\ref{eq:subsetLARC}) contains the vectors $\xi$ and 
$H\xi$. It follows that the LARC condition holds for 
all vectors $\xi$ which are not eigenvectors of $H$.

Clearly $\T$ is not irreducible because any eigenprojection of 
$H$ is an harmonic projection for $\T$.
}
\end{example}


\begin{theorem}\label{th:LARC}
If the LARC holds for some minimal GKSL, then the QMS generated 
by $\Ll$ is irreducible.
\end{theorem} 

\noindent{\bf Proof.} Indeed, if the QMS generated by $\Ll$ is 
not irreducible, for any GKSL representation of $\Ll$ by means 
of operators $G$, $L_\ell$, by Theorem III.1 of \cite{FFRR-jmp02},  
there exists a non-trivial subspace $\mathsf{V}$ invariant for $G$ 
and all $L_\ell$. 
This subspace is also invariant for the operators $\widetilde{G}$ 
and $L_\ell$ because of (\ref{eq:widetildeG}). It follows that, 
for all $\xi\in\mathsf{V}$, the linear manifold generated by 
vectors (\ref{eq:subsetLARC}) is contained in $\mathsf{V}$. 
\hfill $\square$

\medskip
It is worth noticing here that vectors (\ref{eq:subsetLARC}) may 
\emph{not} be contained in  $S(\xi)$. Indeed, there is no reason 
why $\widetilde{G}\xi$ should be contained in $S(\xi)$.

The converse of Theorem \ref{th:LARC} 
is not true, indeed, there exist irreducible QMSs 
with a given GKSL of their generator which do not satisfy the 
LARC condition as shows the following example.

\begin{example}\label{ex:3x3}{\rm
Let $\T$ be the QMS on ${\mathcal{B}}({\mathbb{C}}^{3})$ 
generated by 
\[
{\mathcal{L}}(x) = G^*x + L^* x L + x G 
\]
where $L, H$ are the $3\times 3$ matrices 
\begin{equation}
L = \left[\begin{array}{ccc} 0 & 1 & 0 \\ 
                             -1 & 0 & 0 \\
                             0 & 0 & 0 
 \end{array}\right],  \qquad
H =  \left[\begin{array}{ccc} 0 & 0 & -\mi \\ 
                             0 & 0 & 0 \\
                             \mi & 0 & 0 
\end{array}\right]
\end{equation} 
and 
\[
G=-\frac{1}{2}L^*L - \mi H = 
                       \left[\begin{array}{ccc} -1/2 & 0 & -1 \\ 
                             0 & -1/2 & 0 \\
                            1 & 0 & 0 
\end{array}\right]
\]
Since $L$ is anti-self-adjoint its invariant subspaces are generated 
by eigenvectors $(1,\mi,0), (-1,\mi,0), (0,0,1)$ of $-\mi L$. One 
immediately checks that no one-dimensional or two-dimensional 
subspace generated by these vectors, which is obviously $L^2$ 
invariant, is $H$ invariant, therefore it is not $G$ invariant and 
the QMS $\T$ is irreducible.

Clearly 
\[
\widetilde{G}=G-\frac{1}{2}L^2 = -\mi H 
       = \left[\begin{array}{ccc} 0 & 0 & -1 \\ 
                             0 & 0 & 0 \\
                            1 & 0 & 0 
\end{array}\right]
\]
and $[\,\widetilde{G},L\,]=[\,G,L\,]$. 
Straightforward computations yield 
\[
[\,\widetilde{G},L\,]=  \left[\begin{array}{ccc} 0 & 0 & 0 \\ 
                             0 & 0 & -1 \\
                             0 & 1 & 0 
\end{array}\right],
\]
so that, defining 
\[
X_1 =  [\,\widetilde{G},L\,] \qquad 
X_2 =-\widetilde{G}, \qquad X_3 = -L
\]
we find a basis of the Lie algebra of the rotation group $SO(3)$,
which satisfies the commutation relations
$$
[\, X_1 , X_2 \,] = X_3,  \qquad
[\, X_2 , X_3 \,] = X_1 ,  \qquad
[\, X_3 , X_1 \,] = X_2 .
$$
Consequently iterated commutators of $X_1,X_2, X_3$ do 
not give other operators.
Let $e_1=(1,0,0), e_2=(0,1,0), e_3=(0,0,1)$ be the 
canonical orthonormal basis of $\mathbb{C}^3$. 
It is now immediate to check that 
\begin{eqnarray*}
& & X_1 e_2 = e_3, \quad  X_1 e_3 = -e_2, \quad X_2 e_1 = -e_3 
\\ & &
X_2 e_3 =  e_1, \quad\,  X_3\, e_1 = e_2, \quad\, X_3\, e_2 =-e_1  
\end{eqnarray*}
and $X_k e_k = 0$ for $k=1,2,3$. It follows that the linear 
manifold $S(e_k)$ is two-dimensional for all $k=1,2,3$ and the 
(\ref{eq:subsetLARC}) condition does not hold.
}
\end{example}

\section{Generic QMSs}\label{sect:Generic-QMSs}

In this section we show that irreducibility of  generic QMSs is 
equivalent to the LARC condition with respect to the natural 
GKSL representation of $\Ll$.

Generic QMSs arise in the stochastic limit of a open discrete 
quantum system with generic Hamiltonian, interacting with 
Gaussian fields through a dipole type  interaction (see 
\cite{AcFaHa,CaFaHa,CaSaUm-generic}).  Here, as in the 
previous sections, 
the system space is finite-dimensional $\h=\mathbb{C}^d$ 
with orthonormal basis $(e_k)_{k\le j\le d}$. The operators 
$L_\ell$, in this case labeled by a double index $(\ell,k)$ with 
$\ell\not=k$, are 
\[
L_{\ell k} = \gamma_{\ell k}^{1/2}\ketbra{e_k}{e_\ell}
\]
where are $\gamma_{\ell k}\ge 0$  positive constants and the 
Hamiltonian $H$ is a self-adjoint operator diagonal in the given basis 
whose explicit form is not needed here. The generator $\Ll$ is 
\begin{equation}\label{eq:GeneratGeneric} 
\Ll(x) =  \mi [H,x]+\frac{1}{2} \sum_{\ell\not=k} 
\left(-L_{\ell k}^*L_{\ell k} x + 2 L_{\ell k}^* x L_{\ell k}
- x L_{\ell k}^*L_{\ell k}\right). 
\end{equation}
%

The converse of Theorem \ref{th:LARC} holds for generic QMS.  
\begin{theorem}
A generic QMS is irreducible if and only if the LARC condition holds.
\end{theorem}

\noindent{\bf Proof.} 
The restriction of $\Ll$ to the algebra of diagonal matrices coincides 
with the generator of a time continuous classical Markov chain 
with states ${1,\dots, d}$ and jump rates $\gamma_{\ell k}$ 
(see \cite{AcFaHa,CaFaHa}). 
It is easy to see as in \cite{FaHa} that the QMS generated by (\ref{eq:GeneratGeneric}) is irreducible if and only if the classical 
Markov chain is irreducible, i.e. for all pair of states $\ell,m$ with 
$\ell\not=k$, there exist $n\ge 1$ and states $j_1,\dots, j_n$ 
such that 
\begin{equation}\label{eq:class-irred}
\gamma_{\ell j_1}\gamma_{j_1 j_2}\dots \gamma_{j_n k}>0.
\end{equation}
There is no loss of generality in assuming that each $j_i$ is 
not equal to any $j_1,\dots,j_{i-1}$, together with $j_i \neq \ell$.
Indeed, if $j_i=j_{i^\prime}$, we can 
delete all states $j_{i^\prime},\dots, j_{i-1}$ in the sequence 
of transitions $\ell\to j_1 \to\dots \to j_{n}\to k$. 
Consequently, we find the commutation relation
\[
[L_{j_i\,  j_{i+1}}, L_{j_{i-1}\, j_i}] 
= L_{j_i\,  j_{i+1}}  L_{j_{i-1}\, j_i} 
= \left( \gamma_{j_{i-1}j_i}\gamma_{j_i j_{i+1}}\right)^{1/2}
 \ketbra{e_{j_{i+1}}}{e_{j_{i-1}}}
\]
and compute the iterated commutator
\[
[L_{j_n k},[L_{j_{n-1}j_{n}},\dots ]]
=\left( \gamma_{\ell j_1}\gamma_{j_1 j_2}\dots \gamma_{j_n k}
\right)^{1/2} \ketbra{e_k}{e_\ell}.
\]
For all non-zero $\xi=\sum_{1\le j\le d}\xi_j e_j\in\mathbb{C}^d$,  
choose an $\ell$ such that $\xi_\ell\not=0$ and note that 
\[
[L_{j_n k},[L_{j_{n-1}j_{n}},\dots ]] \xi 
=\xi_\ell \left( \gamma_{\ell j_1}\gamma_{j_1 j_2}
\dots \gamma_{j_n k}\right)^{1/2}e_k.
\] 
It follows that, if the QMS is irreducible, i.e. (\ref{eq:class-irred}) 
holds, then also the LARC condition holds. This, together with 
Theorem \ref{th:LARC}, implies that irreducibility and the 
LARC condition are equivalent for generic QMSs.  
\hfill $\square$


\section*{Acknowledgement}  
The financial support from FONDECYT Grants 1110787 and 1140411 
 is gratefully acknowledged. 

This work was reported by FF at the workshop  ``Operator 
\& Spectral theory, Operator algebras, Non-commutative 
Geometry \& Probability'' in honour of Kalyan B. Sinha - 70th 
Birthday, Kozhikode, India, February 2014. 

The authors would like to take this opportunity 
of wishing a happy birthday to Kalyan thanking him for his 
outstanding scientific contributions, managerial skills  and 
tireless work in dissemination of quantum probability. 
FF is very grateful for the privilege of collaborating with Kalyan on 
quantum stochastic calculus at ISI New Delhi in the early nineties.

FF also would like to express his gratitude for the hospitality 
of the Kerala School of Mathematics and thank Professor 
A.K. Vijayarajan for the enyojable stay  and the very interesting workshop.

\end{document}